\documentclass{gtart}

\def\ifplaintex{\expandafter\ifx\csname documentclass\endcsname\relax}

\def\gtp{{\mathsurround=0pt\it $\cal G\mskip-2mu$eometry \&\ 
$\cal T\!\!$opology $\cal P\!$ublications}}  

\def\recd{{\small Received:\qua\receiveddate\ifx\reviseddate\relax
\else\qquad Revised:\qua\reviseddate\fi\par}} 

\def\lognumber#1{\def\thelognumber{#1}}
\def\volumenumber#1{\def\thevolumenumber{#1}}
\def\volumeyear#1{\def\thevolumeyear{#1}}
\def\papernumber#1{\def\thepapernumber{#1}}
\def\pagenumbers#1#2{\def\startpage{#1}\def\finishpage{#2}}
\def\published#1{\def\publishdate{#1}}

\def\received#1{\def\receiveddate{#1}}
\def\revised#1{\def\reviseddate{#1}}
\def\accepted#1{\def\accepteddate{#1}}

\long\def\asciiabstract#1{\long\def\theasciiabstract{#1}}

\let\\\par\let\thelognumber\relax\let\thevolumenumber\relax
\let\thepapernumber\relax\let\thevolumeyear\relax\let\startpage\relax
\let\finishpage\relax\let\publishdate\relax\let\receiveddate\relax
\let\reviseddate\relax\let\accepteddate\relax\let\theasciititle\relax
\let\theasciiauthors\relax
\let\theasciiabstract\relax

\let\theasciiemail\relax

\ifplaintex
\font\logobig=cmssbx10 scaled 3836
\font\logomed=cmssbx10 scaled 2557
\else
\font\logobig=cmssbx10 scaled 4200
\font\logomed=cmssbx10 scaled 2800
\fi

\long\def\makeagttitle{   
\count0=\startpage
\agt\hfill      
\hbox to 45truept{\vbox to 0pt{\vglue -13truept{\logomed A\kern -.37em{\logobig 
T}\kern -.38em G}\vss}\hss}
\break
{\small Volume \thevolumenumber\ (\thevolumeyear)
\startpage--\finishpage\nl
Published: \publishdate}

\vglue .25truein

{\parskip=0pt\leftskip 0pt plus
1fil\def\\{\par\smallskip}{\Large\bf\thetitle}\par\medskip} \vglue
0.05truein

{\parskip=0pt\leftskip 0pt plus 1fil\def\\{\par}{\sc\theauthors}
\par\medskip}%
 
\vglue 0.03truein 

{\small\leftskip 25truept\rightskip 25truept{\bf Abstract}\stdspace\theabstract

{\bf AMS Classification}\stdspace\theprimaryclass
\ifx\thesecondaryclass\relax\else; \thesecondaryclass\fi\par
{\bf Keywords}\stdspace \thekeywords\par}\vglue 7truept

}   

\ifplaintex
\font\phead=cmsl9 scaled 950
\font\pnum=cmbx10 scaled 913
\font\pfoot=cmsl9 scaled 950
\headline{\vbox to 0pt{\vskip -4.5mm\line{\small\phead\ifnum
\count0=\startpage ISSN 1472-2739 (on-line) 1472-2747 (printed)
\hfill {\pnum\folio}\else\ifodd\count0\def\\{ }%
\ifx\theshorttitle\relax\thetitle\else\theshorttitle\fi\hfill{\pnum\folio}
\else\def\\{ and }{\pnum\folio}\hfill\ifx\theshortauthors\relax\theauthors
\else\theshortauthors\fi\fi\fi}\vss}}
\footline{\vbox to 0pt{\vglue 0mm\line{\small\pfoot\ifnum\count0=\startpage
\copyright\ \gtp\hfill\else
\agt, Volume \thevolumenumber\ (\thevolumeyear)\hfill\fi}\vss}}
\else
\font=cmsl9 scaled 1050
\font\lnum=cmbx10 
\font=cmsl9 scaled 1050
\makeatletter
\def\@oddhead{{\small\lhead\ifnum\count0=\startpage ISSN 1472-2739 
(on-line) 1472-2747 (printed)\hfill {\lnum\number\count0}\else\ifodd\count0
\def\\{ }\ifx\theshorttitle\relax \thetitle \else\theshorttitle\fi\hfill
{\lnum\number\count0}\else\def\\{ and }{\lnum\number\count0}
\hfill\ifx\theshortauthors\relax 
\theauthors\else\theshortauthors\fi\fi\fi}}\def\@evenhead{\@oddhead}
\def\@oddfoot{\small\lfoot\ifnum\count0=\startpage\copyright\ \gtp\hfill\else
\agt, Volume \thevolumenumber\ (\thevolumeyear)\hfill\fi}
\def\@evenfoot{\@oddfoot}
\makeatother
\fi
\let\maketitlepage\makeagttitle
\let\makeshorttitle\maketitlepage
\let\maketitle\maketitlepage

\newwrite\gtoutfile
\long\gdef\makeheadfile{  
{\def\\{, }\def\s{ }
\immediate\openout\gtoutfile head.xxx
\immediate\write\gtoutfile{To: math@arxiv.org}
\immediate\write\gtoutfile{Subject: put OR rep NNNNN:ppppp}
\immediate\write\gtoutfile{--text follows this line--}
\immediate\write\gtoutfile{Proxy-for: \ifx\theasciiauthors\relax
\theauthors\else\theasciiauthors\fi\s<\ifx\theasciiemail\relax\theemail\else\theasciiemail\fi>}
\immediate\write\gtoutfile{\noexpand\\}
\immediate\write\gtoutfile{Authors: \ifx\theasciiauthors\relax
\theauthors\else\theasciiauthors\fi}
{\def\\{ }\immediate\write\gtoutfile{Title: \ifx\theasciititle\relax
\thetitle\else\theasciititle\fi}}
\immediate\write\gtoutfile{Subj-class: GT or SG, GR etc}
\immediate\write\gtoutfile{MSC-class: \theprimaryclass\ifx\thesecondaryclass\relax\else, \thesecondaryclass\fi}
\immediate\write\gtoutfile{Journal-ref: Algebr. Geom. Topol. \thevolumenumber\s
(\thevolumeyear) \startpage-\finishpage}
\immediate\write\gtoutfile{Comments: Published by Algebraic and
Geometric Topology at}
\immediate\write\gtoutfile{\s\s\s  http://www.maths.warwick.ac.uk/agt/AGTVol\thevolumenumber/agt-\thevolumenumber-\thepapernumber.abs.html}
\immediate\write\gtoutfile{\noexpand\\}
\immediate\write\gtoutfile{}
\ifx\theasciiabstract\relax
\immediate\write\gtoutfile{\theabstract}\else
\immediate\write\gtoutfile{\theasciiabstract}\fi
\immediate\write\gtoutfile{}
\immediate\write\gtoutfile{\noexpand\\}
\immediate\write\gtoutfile{}
\immediate\closeout\gtoutfile}}  

\def\maketitlepage{\makeagttitle\makeheadfile}
\let\makeshorttitle\maketitlepage
\let\maketitle\maketitlepage

\def\ifplaintex{\expandafter\ifx\csname documentclass\endcsname\relax}

\def\gtp{{\mathsurround=0pt\it $\cal G\mskip-2mu$eometry \&\ 
$\cal T\!\!$opology $\cal P\!$ublications}}  

\def\recd{{\small Received:\qua\receiveddate\ifx\reviseddate\relax
\else\qquad Revised:\qua\reviseddate\fi\par}} 

\def\lognumber#1{\def\thelognumber{#1}}
\def\volumenumber#1{\def\thevolumenumber{#1}}
\def\volumeyear#1{\def\thevolumeyear{#1}}
\def\papernumber#1{\def\thepapernumber{#1}}
\def\pagenumbers#1#2{\def\startpage{#1}\def\finishpage{#2}}
\def\published#1{\def\publishdate{#1}}

\def\received#1{\def\receiveddate{#1}}
\def\revised#1{\def\reviseddate{#1}}
\def\accepted#1{\def\accepteddate{#1}}

\long\def\asciiabstract#1{\long\def\theasciiabstract{#1}}

\let\\\par\let\thelognumber\relax\let\thevolumenumber\relax
\let\thepapernumber\relax\let\thevolumeyear\relax\let\startpage\relax
\let\finishpage\relax\let\publishdate\relax\let\receiveddate\relax
\let\reviseddate\relax\let\accepteddate\relax\let\theasciititle\relax
\let\theasciiauthors\relax
\let\theasciiabstract\relax

\let\theasciiemail\relax

\ifplaintex
\font\logobig=cmssbx10 scaled 3836
\font\logomed=cmssbx10 scaled 2557
\else
\font\logobig=cmssbx10 scaled 4200
\font\logomed=cmssbx10 scaled 2800
\fi

\long\def\makeagttitle{   
\count0=\startpage
\agt\hfill      
\hbox to 45truept{\vbox to 0pt{\vglue -13truept{\logomed A\kern -.37em{\logobig 
T}\kern -.38em G}\vss}\hss}
\break
{\small Volume \thevolumenumber\ (\thevolumeyear)
\startpage--\finishpage\nl
Published: \publishdate}

\vglue .25truein

{\parskip=0pt\leftskip 0pt plus
1fil\def\\{\par\smallskip}{\Large\bf\thetitle}\par\medskip} \vglue
0.05truein

{\parskip=0pt\leftskip 0pt plus 1fil\def\\{\par}{\sc\theauthors}
\par\medskip}%
 
\vglue 0.03truein 

{\small\leftskip 25truept\rightskip 25truept{\bf Abstract}\stdspace\theabstract

{\bf AMS Classification}\stdspace\theprimaryclass
\ifx\thesecondaryclass\relax\else; \thesecondaryclass\fi\par
{\bf Keywords}\stdspace \thekeywords\par}\vglue 7truept

}   

\ifplaintex
\font\phead=cmsl9 scaled 950
\font\pnum=cmbx10 scaled 913
\font\pfoot=cmsl9 scaled 950
\headline{\vbox to 0pt{\vskip -4.5mm\line{\small\phead\ifnum
\count0=\startpage ISSN 1472-2739 (on-line) 1472-2747 (printed)
\hfill {\pnum\folio}\else\ifodd\count0\def\\{ }%
\ifx\theshorttitle\relax\thetitle\else\theshorttitle\fi\hfill{\pnum\folio}
\else\def\\{ and }{\pnum\folio}\hfill\ifx\theshortauthors\relax\theauthors
\else\theshortauthors\fi\fi\fi}\vss}}
\footline{\vbox to 0pt{\vglue 0mm\line{\small\pfoot\ifnum\count0=\startpage
\copyright\ \gtp\hfill\else
\agt, Volume \thevolumenumber\ (\thevolumeyear)\hfill\fi}\vss}}
\else
\font=cmsl9 scaled 1050
\font\lnum=cmbx10 
\font=cmsl9 scaled 1050
\makeatletter
\def\@oddhead{{\small\lhead\ifnum\count0=\startpage ISSN 1472-2739 
(on-line) 1472-2747 (printed)\hfill {\lnum\number\count0}\else\ifodd\count0
\def\\{ }\ifx\theshorttitle\relax \thetitle \else\theshorttitle\fi\hfill
{\lnum\number\count0}\else\def\\{ and }{\lnum\number\count0}
\hfill\ifx\theshortauthors\relax 
\theauthors\else\theshortauthors\fi\fi\fi}}\def\@evenhead{\@oddhead}
\def\@oddfoot{\small\lfoot\ifnum\count0=\startpage\copyright\ \gtp\hfill\else
\agt, Volume \thevolumenumber\ (\thevolumeyear)\hfill\fi}
\def\@evenfoot{\@oddfoot}
\makeatother
\fi
\let\maketitlepage\makeagttitle
\let\makeshorttitle\maketitlepage
\let\maketitle\maketitlepage

\newwrite\gtoutfile
\long\gdef\makeheadfile{  
{\def\\{, }\def\s{ }
\immediate\openout\gtoutfile head.xxx
\immediate\write\gtoutfile{To: math@arxiv.org}
\immediate\write\gtoutfile{Subject: put OR rep NNNNN:ppppp}
\immediate\write\gtoutfile{--text follows this line--}
\immediate\write\gtoutfile{Proxy-for: \ifx\theasciiauthors\relax
\theauthors\else\theasciiauthors\fi\s<\ifx\theasciiemail\relax\theemail\else\theasciiemail\fi>}
\immediate\write\gtoutfile{\noexpand\\}
\immediate\write\gtoutfile{Authors: \ifx\theasciiauthors\relax
\theauthors\else\theasciiauthors\fi}
{\def\\{ }\immediate\write\gtoutfile{Title: \ifx\theasciititle\relax
\thetitle\else\theasciititle\fi}}
\immediate\write\gtoutfile{Subj-class: GT or SG, GR etc}
\immediate\write\gtoutfile{MSC-class: \theprimaryclass\ifx\thesecondaryclass\relax\else, \thesecondaryclass\fi}
\immediate\write\gtoutfile{Journal-ref: Algebr. Geom. Topol. \thevolumenumber\s
(\thevolumeyear) \startpage-\finishpage}
\immediate\write\gtoutfile{Comments: Published by Algebraic and
Geometric Topology at}
\immediate\write\gtoutfile{\s\s\s  http://www.maths.warwick.ac.uk/agt/AGTVol\thevolumenumber/agt-\thevolumenumber-\thepapernumber.abs.html}
\immediate\write\gtoutfile{\noexpand\\}
\immediate\write\gtoutfile{}
\ifx\theasciiabstract\relax
\immediate\write\gtoutfile{\theabstract}\else
\immediate\write\gtoutfile{\theasciiabstract}\fi
\immediate\write\gtoutfile{}
\immediate\write\gtoutfile{\noexpand\\}
\immediate\write\gtoutfile{}
\immediate\closeout\gtoutfile}}  

\def\maketitlepage{\makeagttitle\makeheadfile}
\let\makeshorttitle\maketitlepage
\let\maketitle\maketitlepage

\lognumber{29}
\volumenumber{1}
\volumeyear{2001}
\papernumber{29}
\pagenumbers{579}{587}
\received{17 May 2001}
\revised{15 August 2001}
\accepted{11 October 2001}
\published{18 October 2001}

\usepackage{amssymb}
\usepackage{amsmath}
\usepackage{amsthm}

\let\Bbb\mathbb

\begin{document}

\newtheorem{thm}{Theorem}[section]
\newtheorem{lem}[thm]{Lemma}
\newtheorem{cor}[thm]{Corollary}
\newtheorem{conj}[thm]{Conjecture}
\newtheorem{cond}{Condition}
\newtheorem{qn}[thm]{Question}

\theoremstyle{definition}
\newtheorem{defn}[thm]{Definition}
\newtheorem{exa}[thm]{Example}

\theoremstyle{remark}
\newtheorem{rmk}[thm]{Remark}

\def\R{\Bbb R}
\def\Z{\Bbb Z}
\def\CP{\Bbb {CP}}
\def\H{\Bbb H}
\def\F{\mathcal F}
\def\C{\Bbb C}
\def\A{\bf A}
\def\homeo{\text{Homeo}}
\def\u{{\text{univ}}}
\def\til{\widetilde}
\def\hat{\widehat}
\def\map{\text{Map}}

\newenvironment{pf}{\proof}{\endproof}

\title{Leafwise smoothing laminations}

\author{Danny Calegari}
\address{Department of Mathematics\\Harvard\\Cambridge, MA 02138}
\email{dannyc@math.harvard.edu}
\url{www.math.harvard.edu/\char'176 dannyc}

\begin{abstract}
We show that every topological surface lamination $\Lambda$ of a $3$--manifold
$M$ is isotopic to one with smoothly immersed leaves. This carries out a
project proposed by Gabai in \cite{dG93}. Consequently any such lamination
admits the structure of a {\em Riemann surface lamination}, and therefore
useful structure theorems of Candel \cite{aC93} and Ghys \cite{eG99} apply.
\end{abstract}
\asciiabstract{ We show that every topological surface lamination of a
3-manifold M is isotopic to one with smoothly immersed leaves. This
carries out a project proposed by Gabai in [Problems in foliations and
laminations, AMS/IP Stud. Adv. Math. 2.2 1--33]. Consequently any such
lamination admits the structure of a Riemann surface lamination, and
therefore useful structure theorems of Candel [Uniformization of
surface laminations, Ann. Sci. Ecole Norm. Sup. 26 (1993) 489--516]
and Ghys [Dynamique et geometrie complexes, Panoramas et Syntheses 8
(1999)] apply.}  \primaryclass{57M50} \keywords{Lamination, foliation,
leafwise smooth, 3--manifold}
\maketitlepage

\section{Basic notions}

\begin{defn}
A {\em lamination} is a topological space which can be covered by open
charts $U_i$ with a local product structure $\phi_i:U_i \to \R^n \times X$ 
in such a way that the manifold--like factor is preserved in the
overlaps. That is, for $U_i \cap U_j$ nonempty,
$$\phi_j \circ \phi_i^{-1}: \R^n \times X \to \R^n \times X$$
is of the form $$\phi_j \circ \phi_i^{-1}(t,x) = (f(t,x),g(x))$$
\end{defn}
The maximal continuations of the local manifold slices $\R^n \times \text{point}$
are the {\em leaves} of the lamination. A {\em surface lamination} is
a lamination locally modeled on $\R^2 \times X$. We usually assume that
$X$ is locally compact.

\begin{defn}
A lamination is {\em leafwise $C^n$} for $n\ge 2$ if the leafwise
transition functions $f(t,x)$ can be chosen in such a way
that the mixed partial derivatives in $t$ of orders less than or equal
to $n$ exist for each $x$, and vary continuously as functions of $x$.
\end{defn}

A {\em leafwise $C^n$ structure} on a lamination $\Lambda$ induces on
each leaf $\lambda$ of $\Lambda$ a $C^n$ manifold structure, in the
usual sense.

\begin{defn}
An embedding of a leafwise $C^n$ lamination $i:\Lambda \to M$ 
into a manifold $M$ is an {\em $C^n$ immersion}
if, for some $C^n$ structure on $M$,
for each leaf $\lambda$ of $\Lambda$ the embedding
$\lambda \to M$ is $C^n$.
\end{defn}

Note that if $i:\Lambda \to M$ is an embedding with the property that the
image of each leaf $i(\lambda)$ is locally a $C^n$ submanifold, and these
local submanifolds vary continuously in the $C^n$ topology, then
there is a unique leafwise $C^n$ structure on $\Lambda$ for which $i$
is a $C^n$ immersion.

A foliation of a manifold is an example of a lamination. For a foliation
to be leafwise $C^n$ is {\it a priori} weaker than to ask for it to be
$C^n$ immersed.

\begin{exa}
Let $M$ be a manifold which is not stably smoothable, and $N$ a compact
smooth manifold.
Then $M\times N$ has the structure of a leafwise smooth
foliation (by parallel copies of $N$), but there is no smooth structure
on $M \times N$ for which the embedding of the foliation is a smooth immersion,
since there is no smooth structure on $M \times N$ at all.
\end{exa}

\begin{rmk}
For readers unfamiliar with the notion, the ``tangent bundle'' of a topological
manifold (i.e. a regular neighborhood of the diagonal in $M \times M$) is stably
(in the sense of $K$--theory) classified by a homotopy class of maps
$f:M \to BTOP$ for a certain topological space $BTOP$. There is a fibration
$p:BO \to BTOP$, and the problem of lifting $f$ to $\hat{f}:M \to BO$
such that $p\hat{f} = f$ represents an obstruction to finding a smooth structure
on $M$. For $N$ smooth as above, the composition
$$M \to M \times \text{point} \subset M \times N \to BTOP$$ is homotopic to $f$, and therefore no lifting
of the structure exists on $M\times N$ if none existed on $M$. For a reference,
see \cite{KS77}, or the very readable \cite{yR01}.
\end{rmk}

With notation as above, the tangential quality of $\F$ is controlled by the quality
of $f(\cdot,x)$ for each fixed $x$, for $f$ the first component of a transition
function.
For sufficiently large $k$ and $n-k$
questions of ambiently smoothing {\em foliated manifolds} come down to obstruction
theory and classical surgery theory, as for example in \cite{KS77}.
But in low dimensions, the situation is more elementary and more
hands--on.

\section{Some $3$--manifold topology}

Let $M$ be a topological $3$--manifold. It is a classical theorem of Moise
(see \cite{Moise}) that $M$
admits a PL or smooth structure, unique up to conjugacy.

\begin{lem}\label{generic_family}
Let $\Sigma$ be a topological surface.
Let $S^1_j$ be a countable collection of circles, and let
$\Phi:\coprod_j S^1_j \times I \to \Sigma$ be a map with the following
properties:
\begin{enumerate}
\item{For each $t \in I$, $\Phi(\cdot,t):S^1_j \to \Sigma$ is
an embedding.}
\item{For each $t \in I$ and each pair $j,k$ the intersection
$$\Phi(S^1_j,t) \cap \Phi(S^1_k,t)$$ is finite, and its cardinality is constant
as a function of $t$ away from finitely many values.}
\item{For every compact subset $K \subset \Sigma$ the set of $j$ for which
$\Phi(S^1_j,t) \cap K$ is nonempty for some $t$ is finite.}
\end{enumerate}
Then there is a PL (resp. smooth) structure on $\Sigma \times I$ such that
the graph of each map $\Gamma_j(\Phi):S^1_j \times I \to \Sigma \times I$ is PL (resp. smooth).
\end{lem}
Here the {\em graph $\Gamma_j(\Phi)$ of $\Phi$} is the function
$\Gamma_j(\Phi):S^1_j \times I \to \Sigma \times I$ defined by
$$\Gamma_j(\Phi)(\theta,t) = (\Phi(\theta),t)$$
\begin{pf}
The conditions imply that the image of $\coprod_j S^1_j$ in $\Sigma$ for a fixed
$t$ is topologically a locally finite graph. Such a structure in a $2$ manifold
is locally flat, and the combinatorics of any finite subgraph is locally constant
away from isolated values of $t$. It is therefore straightforward to construct
a PL (resp. smooth) structure on a collar neighborhood of the image of
$\coprod_j S^1_j \times I$ in $\Sigma \times I$. This can be extended canonically
to a PL (resp. smooth) structure on $\Sigma \times I$, by the relative version
of Moise's theorem (see \cite{Moise}).
\end{pf}

\begin{lem}\label{generic_interpolation}
Let $\Phi:S^1_j \times I \to \Sigma$ satisfy the conditions of 
lemma~\ref{generic_family}. Let $\Psi_0:S^1 \to \Sigma$ and $\Psi_1:S^1 \to \Sigma$
be homotopic embeddings such that $\Psi_0(S^1)$ intersects finitely many
circles in $\Phi(\cdot,0)$ in finitely many points, and similarly for $\Psi_1(S^1)$.
Then there is a map $\Psi:S^1 \times I \to \Sigma$ which is a homotopy
between $\Psi_0$ and $\Psi_1$ so that 
$$\Phi\cup\Psi:\Bigl(\coprod_i S^1_i \coprod S^1\Bigr) \times I \to \Sigma$$
satisfies the conditions of lemma~\ref{generic_family}.
\end{lem}
\begin{pf}
Since the combinatorics of the image of $\Phi$ is locally finite,
and since the image of $\Psi$ is bounded,
it suffices to treat the case when $\Phi$ is constant as a function of $t$.

Choose a PL structure on $\Sigma$ for which the image of
$\Phi(\cdot,0)$ and $\Psi_0$ are polygonal. Then produce a polygonal homotopy from
$\Psi_0$ (with respect to this polygonal structure) to a new polygonal $\Psi'_0$
such that $\Psi'_0(S^1)$ and $\Psi_1(S^1)$ intersect the image of
$\Phi(\cdot,t)$ in a finite set of points in the same combinatorial configuration.
Then $\Psi'_0$ is isotopic to $\Psi_1$ rel. its intersection with the image of
$\Phi(\cdot,t)$.
\end{pf}

\section{Surface laminations of $3$--manifolds}

\begin{defn}
Let $\F$ be a codimension one foliation of a $3$--manifold $M$.
A {\em snake} in $M$ is an embedding
$\phi:D^2 \times I \to M$
where $D^2$ denotes the open unit disk, and $I$ the open unit interval,
which extends to an embedding of the closure of $D^2 \times I$, in such
a way that each horizontal disk gets mapped into a leaf $\lambda$ of $\F$.
That is, $\phi:D^2 \times t \to \lambda$.
\end{defn}

The terminology suggests that we are typically interested in snakes which
are reasonably small and thin in the leafwise direction, and possibly large
in the transverse direction.

A collection of snakes in a foliated manifold intersect a leaf $\lambda$
of $\F$ in a locally finite collection of open disks. For a snake $S$,
let $\partial_v \overline{S}$ denote the ``vertical boundary'' of the
closed ball $\overline{S}$; this is topologically an embedded closed cylinder
transverse to $\F$, intersecting each leaf in an inessential circle.

We say that an open cover of $M$ by finitely many snakes $S_i$ is
{\em combinatorially tame} if the embeddings $\partial_v \overline{S_i} \to M$
are locally of the form described in lemma~\ref{generic_family}.

Note that the induced pattern on each leaf $\lambda$ of $\F$ of the circles
$\partial_v \overline{S_i} \cap \lambda$ is topologically conjugate to the
transverse intersection of a locally finite collection of polygons.

\begin{lem}\label{tame_snakes}
A codimension one foliation $\F$ of a closed $3$--manifold $M$ admits a combinatorially
tame open cover by finitely many snakes.
\end{lem}
\begin{pf}
Since $M$ is compact, any cover by snakes contains a finite subcover;
any such cover induces a locally finite cover of each leaf.
We prove the lemma by induction.

Let $S_i$ be a collection of snakes in $M$ which is combinatorially tame.
Let $C_i = \partial_v \overline{S_i}$ be their vertical boundaries,
and let $S$ be another snake with vertical boundary $C$. We will show that there
is a snake $S'$ containing $S$ such that the collection $\lbrace S_i \rbrace \cup
\lbrace S' \rbrace$ is combinatorially tame.

Let $\lambda_t$ for $t\in I$ parameterize the foliation of
$\overline{S}$. Let $E_i(t)$ denote the pattern of
circles $C_i \cap \lambda_t$ in a neighborhood of $E(t) = C\cap \lambda_t$.
By hypothesis, the $C_i$ can be thought of as polygons with respect to a
PL structure on $\lambda_t$. Then $E(t)$ can be {\em straightened} to a polygon
$E(t)'$ in general position with respect to the $E_i(t)$ in a small neighborhood,
where the interior of the region in $\lambda_t$ bounded by $E(t)'$ contains $E(t)$.
If $\lambda_t$ does not intersect the horizontal boundary of any $\overline{S_i}$,
then the combinatorial pattern of intersections of the $E_i(t)$ is
locally generic --- i.e. the pattern might change, but it changes by the
graph of a generic PL isotopy, by lemma~\ref{generic_family}.

It follows that we can extend the straightening of
$E(t)$ to $E(t)'$ for some collar neighborhood of $t=0$. In general,
a straightening of $E(t)$ to $E(t)'$ can be extended in the positive
direction until a $t_0$ which contains some lower horizontal
boundary of an $\overline{S_i}$. The straightening can be extended past an
upper horizontal boundary of an $\overline{S_i}$ without any problems, since
the combinatorial pattern of intersections becomes simpler: circles disappear.

The straightening of $E(t)$ over all $t$ can be done by {\em welding} straightenings
centered at the finitely many values of $t$ which contain horizontal boundary
of some $\overline{S_i}$. Call these critical values $t_j$. So we can
produce a finite collection of straightenings $E(t) \to E(t)'_j$ each
valid on the open interval $t \in (t_{j-1}, t_{j+1})$. To weld these
straightenings together at intermediate values $s_j$ where $t_j < s_j < t_{j+1}$,
we insert a PL isotopy from $E(s_j)'_j$ to $E(s_j)'_{j+1}$ in a little
collar neighborhood of $s_j$, by appealing to lemma~\ref{generic_interpolation}.
So these welded straightenings
give a straightening of $E(t)$ for all $t \in I$, and they bound a snake
$S'$ with the requisite properties.

To prove the lemma, cover $M$ with finitely many snakes $S_i$, and apply the
induction step to straighten $S_j$ while fixing $S_k$ with $k<j$.
Since snakes can be straightened by an arbitrarily small (in the $C^0$ topology)
homotopy, the union of straightened snakes can also be made to cover $M$, and
we are done.
\end{pf}

\begin{lem}\label{smooth_immersion}
Let $M$ be a $3$--manifold, and $\F$ a foliation of $M$ by surfaces. Then $\F$
is isotopic to a foliation such that all leaves are PL or smoothly immersed,
and the images of leaves vary locally continuously in the $C^\infty$ topology.
\end{lem}
\begin{pf}
If $S_i$ is a combinatorially tame cover of $\F$ by snakes, the image of
the union $\cup_i \partial \overline{S_i}$ can be taken to be a PL or smooth
2 complex $\Sigma$ 
in $M$, whose complementary regions are polyhedral 3 manifolds.
Each complementary region is foliated as a product by $\F$. We can straighten
$\F$ cell--wise inductively on its intersection with the skeleta of $\Sigma$.
First, we keep $\F \cap \Sigma^1$ constant. Then the foliation of
$\F \cap (\Sigma^2\backslash \Sigma^1)$
by lines can be straightened to be PL or smooth, and
this straightened foliation extended in a PL or smooth manner over the
product complementary regions in $M - \Sigma$.
\end{pf}

\begin{thm}\label{lamination_smooth_immersion}
Let $\Lambda$ be a surface lamination in a $3$--manifold $M$. Then $\Lambda$ is
isotopic to a lamination such that all leaves are PL or smoothly immersed, and
the images of leaves vary locally continuously in the $C^\infty$ topology.
\end{thm}
\begin{pf}
By the definition of a lamination, there is an open cover of
$\Lambda$ by balls $B_i$
such that $\Lambda \cap B_i$ is a product lamination, which can be extended to a
product foliation. It is straightforward to produce an open submanifold $N$ with
$\Lambda \subset N \subset M$ such that $N$ can be foliated by a foliation $\F$
which contains $\Lambda$ as a sublamination. Then the open manifold $N$ can be
given a PL or smooth structure in which $\F$, and hence $\Lambda$, is PL or
smoothly immersed, by lemma~\ref{smooth_immersion}. This PL or smooth structure can
be extended compatibly over $M-N$ by Moise's theorem.
\end{pf}

\begin{cor}
Let $\Lambda$ be a surface lamination in a $3$--manifold $M$. Then $\Lambda$ admits
a leafwise PL or smooth structure.
\end{cor}

In particular, such a lamination admits the structure of a Riemannian surface
lamination. In Gabai's problem list \cite{dG93}, he lists 
theorem~\ref{lamination_smooth_immersion}
as a ``project''. The corollary allows us to apply the technology of complex analysis and
algebraic geometry to such laminations; in particular, the following theorems
of Candel and Ghys from \cite{aC93} and \cite{eG99} apply:

\begin{thm}[Candel]
Let $\F$ be an essential Riemann surface lamination of an atoroidal
$3$--manifold. Then there
exists a continuously varying path metric on $\F$ for which the leaves of $\F$
are locally isometric to $\H^2$.
\end{thm}

\begin{thm}[Ghys]
Let $\F$ be a taut foliation of a $3$--manifold $M$ with Riemann surface leaves.
Then there is an embedding $e:M \to \CP^n$ for some $n$ which is leafwise
holomorphic. That is, $e:\lambda \to \CP^n$ is holomorphic for each leaf $\lambda$.
\end{thm}

\vfill
\pagebreak

\Addresses
\recd

\end{document}